\documentclass[aps,pre,preprint,groupedaddress,showpacs]{revtex4}

\usepackage{graphicx}
\usepackage{dcolumn}
\usepackage{bm}
\usepackage{subfigure}

\begin{document}


\title{Evidence of Long Range Order in the Riemann Zeta Function}


\author{Ronald Fisch}
\email[]{ronf124@yahoo.com}
\affiliation{382 Willowbrook Dr.\\
North Brunswick, NJ 08902}


\date{\today}

\begin{abstract}
We have done a statistical analysis of some properties of the contour
lines Im$(\zeta ( s ))$ = 0 of the Riemann zeta function.  We find that
this function is broken up into strips whose average width on the
critical line does not appear to vary with height.  We also compute the
position of the primary zero for the lowest 200 strips, and find that
this probability distribution also appears to be scale invariant.

\end{abstract}


\maketitle

\section{Introduction}

Let $s ~=~ \sigma ~+~ it$, with $\sigma$ and $t$ real variables.
Then for $\sigma ~>~ 1$ the Riemann zeta function is defined as
\begin{equation}
  \zeta ( s ) ~=~  \sum_{n=1}^{\infty} n^{-s}  \, .
\end{equation}

It follows immediately from Eqn. (1), that for any $t$
\begin{equation}
  \lim_{\sigma \to +\infty} \zeta ( s ) ~=~ 1  \, .
\end{equation}
It was shown by Riemann\cite{BCRW08} that $\zeta ( s )$ can be
analytically continued to a function which is meromorphic in the
complex plane, that its only divergence is a simple pole at $s$ = 1,
and that it has no zeroes on the half-plane $\sigma ~>~ 1$.

The Riemann Hypothesis\cite{BCRW08} (RH) is that the only zeroes of
$\zeta ( s )$ which do not lie on the real axis lie on the critical
line $s ~=~ {1/2} ~+~ it$.  It has resisted rigorous proof for over
150 years, and is now widely considered to be the most important
unsolved problem in mathematics.\cite{Con03}  The significance of
RH for physics has been shown by many authors.\cite{BK99,ST08,
MDMSWS10,SH11,FHK12}

In this work we will study the statistical properties of the contour
lines Im$(\zeta ( s ))$ = 0.  We will find that there appears to be
a type of long range order in the spacing of these contour lines.
We will also find a statistical quantity whose second moment appears
to be independent of the height $t$.  In other words, this quantity
has a scaling exponent of zero, which strongly suggests the existence
of some hidden conservation law in $\zeta ( s )$.

\section{Numerical results}

The contour lines Im$(\zeta ( s ))$ = 0 and Re$(\zeta ( s ))$ = 0 are
highly constrained by the Cauchy-Riemann equations.  An extensive
discussion of their properties, including many illustrations, has been
given by Arias de Reyna.\cite{Arias03}  For this work, the author used
the computer program of Collins.\cite{Col09}  Due to Eqn. (2) and the
Cauchy-Riemann equations, there are contour lines Im$(\zeta ( s ))$ = 0
which break up $\zeta ( s )$ into what appear to be essentially horizontal
strips, although these contours are not actually straight lines.

If we define the phase $\theta ( s )$ by
\begin{equation}
  \zeta ( s ) ~=~ |\zeta ( s )| \exp ( i \theta ( s ))  \, ,
\end{equation}
then the contour lines at the top and the bottom of each strip have
$\cos(\theta)$ = 1, {\it i.e.} $\theta$ = 0.  Reading the crossings of
the critical line, $s$ = 1/2, by these contour lines, we plot the number
of strips as a function of the height $t$ in Fig.~1 for the first 200
strips.  Because the crossings were measured by eye, the heights were
rounded off to integer values.  This is merely a rounding error.  It
does not accumulate.  Therefore this rounding error has virtually no
effect on the slope of the least squares fitting line shown in Fig.~1.

\begin{figure}
\includegraphics[width=3.4in]{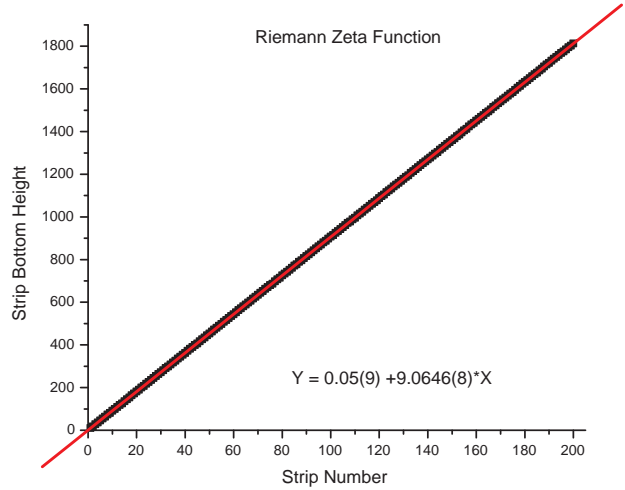}
\caption{\label{Fig.1}  Height of the bottom of a strip (rounded to the
nearest integer) on the critical line $\sigma$ = 1/2 as a function of
strip number for the first 200 strips.}
\end{figure}

The linear least squares fit to the data is
\begin{equation}
  Y ~=~ 0.05(9) ~+~ 9.0646(8)*X  \, ,
\end{equation}
where the numbers in parentheses are the statistical errors in the last
significant figure.  Since there is a nontrivial distribution of strip
widths, there is a small amount of jitter of the data around the fitting
line.  However, there is absolutely no indication of any curvature in
the fit.  The bottom of the first strip sits at a height of approximately
$t$ = 10.  It is therefore somewhat mysterious that the Y-intercept of
the fitting line is consistent (within the statistical error) with a
value of zero.  If one fits the heights of the tops of the strips instead,
one finds (unsurprisingly) that the slope of the fitting line is
essentially unchanged, but the Y-intercept is now found to be 9.07(9).

For large positive $\sigma$ Eqn. (1) can be approximated by its first two
terms.  Under this condition
\begin{equation}
  {\rm Im} (\zeta ( s )) ~\approx~  {\rm Im} (2^{-s}) ~=~ 2^{-\sigma} \sin (\ln(2)t)  \, .
\end{equation}
Thus the strip boundaries for large positive $\sigma$ and $t > 0$ will be
\begin{equation}
  t ~\approx~ 2 m \pi /\ln(2)  \, ,
\end{equation}
where the strip number, $m$, is a positive integer.  The numerical value of
$2 \pi/\ln(2)$ is 9.06472... .  Assuming the RH is correct, it seems a
reasonable conjecture that the strips remain essentially horizontal for any
value of $t$ when $\sigma > 0$, which implies that the slope of the least
squares fit for the heights of the bottom of each strip will be independent
of $\sigma$.  The author sees no reason, however, why the Y-intercept of this
fit should be independent of $\sigma$.  In fact, it appears that for
$\sigma < 0$ this Y-intercept becomes clearly greater than zero.

The reader should note that, since the sum on the right hand side of Eqn. (1)
does not even converge for $\sigma$ = 1/2, it is quite remarkable for the
data taken on the critical line, shown in Fig.~1, to be well fit by a straight
line with a slope of $2 \pi/\ln(2)$.  Based on the analysis of Berry and
Keating,\cite{BK99} for example, one might have expected to see oscillations
about this line coming from the higher order terms of the sum.

\begin{figure}
\includegraphics[width=3.4in]{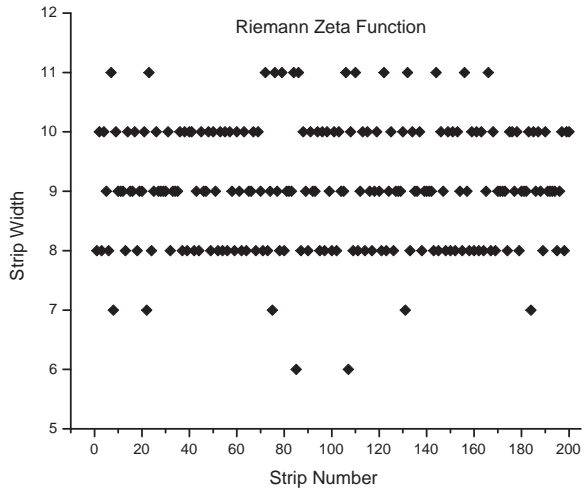}
\caption{\label{Fig.2}  Strip width (rounded to integer values) on the critical
line $\sigma$ = 1/2 as a function of strip number, for the first 200 strips.}
\end{figure}

The strip widths, measured at $\sigma$ = 1/2, as a function of strip number
are shown in Fig.~2.  As mentioned above, in the process of measuring by
eye, these widths were rounded off to integer values.  The distribution of
strip widths appears to be independent of height.  Due to the measurement
process, we cannot determine how well this distribution would be fit by a
Poisson distribution.  There is a clear tendency for the distribution to
become narrower as $\sigma$ becomes greater than 1, i.e. the narrow strips
become wider and the wide strips become narrower.  It would be interesting to
calculate the analytical form of this distribution for large positive $\sigma$.

The reader should note that Arias de Reyna\cite{Arias03} displays a few strips
near the critical line at much greater values of $t$ than those studied here,
and the widths of those strips are consistent with our estimates of the widths
for $t ~<~ 2000$.

\begin{figure}
\includegraphics[width=3.4in]{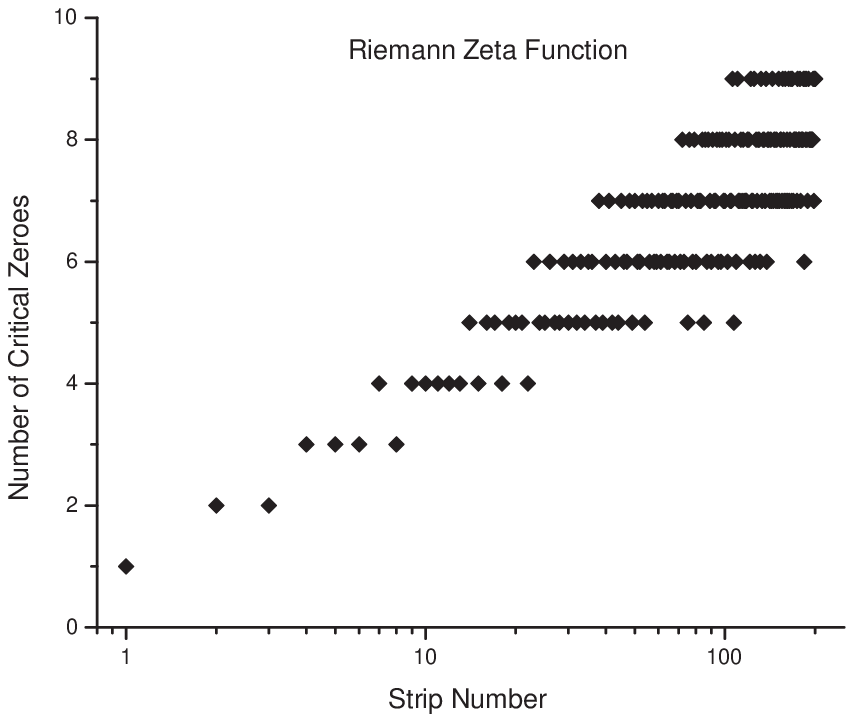}
\caption{\label{Fig.3}  Number of critical zeroes in a strip versus strip
number, for the first 200 strips. The $X$-axis is scaled logarithmically.}
\end{figure}

The number of critical zeroes in a strip versus strip height is shown in
Fig.~3.  We see that the number of zeroes increases approximately
logarithmically with strip number.  Because we saw in Fig.~1 and Fig.~2
that the distribution of strip widths appears to be independent of the
strip number, this is unsurprising.  It has been known for many years
that the density of the critical zeroes increases approximately as
$\log ( t )$.  Actually, the average density of critical zeroes as a
function of $t$ is known to a much greater precision than this.\cite{Arias03}

\begin{figure}
\includegraphics[width=3.4in]{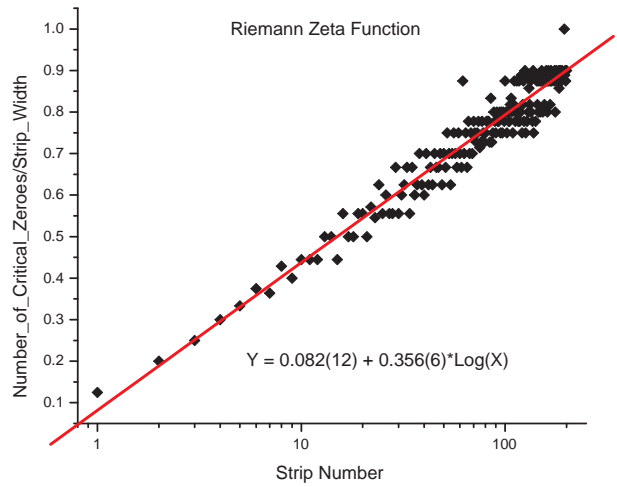}
\caption{\label{Fig.4}  (Number of zeroes on the strip)/(strip width)
versus strip number, for the first 200 strips. The $X$-axis is scaled
logarithmically.}
\end{figure}

The width of each strip on the line $\sigma$ = 1/2 is highly correlated
with the number of zeroes in the strip.  To see this, in Fig.~4 we display
the function (number of zeroes on the strip) divided by (strip width) versus
the strip number.  Comparing Fig.~4 to Fig.~3, it is clear that the scatter
in Fig.~4 is much smaller.  This is true despite the large roundoff error
in our measurement of the strip width.

\begin{figure}
\includegraphics[width=3.4in]{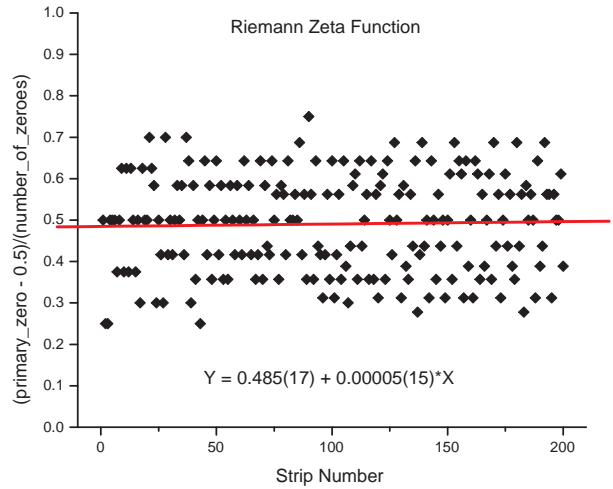}
\caption{\label{Fig.5}  (Number of the primary zero $-~0.5$)/(Number of
zeroes) versus strip number, for the first 200 strips.}
\end{figure}

Due to the requirements of the Cauchy-Riemann equations, in each strip
there is one special zero, which we will call the primary zero.  Each
primary zero has the property that the contour with phase $\theta$ = 0
going out of it extends to $\sigma ~=~ +\infty$, at a height
\begin{equation}
  t ~\approx~ (2 m + 1) \pi /\ln(2)  \, .
\end{equation}
For all the other critical zeroes, the contour with phase $\theta ~=~ 0$
goes to $\sigma ~=~ -\infty$.

One can now ask the question ``Where is the primary zero located in the
strip?"  It seems obvious, by reason of symmetry, that the average
position of the primary zero should be at the center of the strip.
However, there is no symmetry reason why the probability distribution
for the primary zero should be uniform.  In Fig.~5 we display the values
for the function (number of the primary zero (counting from the bottom
of the strip) $-~ 0.5$) divided by (number of zeroes in the strip) versus
the strip number.  The subtraction of 0.5 in the numerator is necessary
so that this function has the value 0.5 when the primary zero is the
middle zero.

The linear least squares fit to the data of Fig.~5 shows that the average
position of the primary zero is indeed at the center of the strip.
Remarkably, one sees that the value of the variance of this probability
distribution seems to be independent of the strip height.  This
observation is confirmed by comparing the variance for the first
100 strips with the variance for the second 100 strips.

If this probability distribution remains nontrivial ({\it i.e.} neither
becoming uniform nor collapsing) in the limit $t \to \infty$, then we must
conclude that all of the zeroes in a strip are highly correlated with each
other.  Since the number of zeroes in a strip is expected to diverge as
$t \to \infty$, this is a very interesting effect.  It strongly suggests
that there must exist some unknown (to the author) symmetry or hidden
conservation law at work.

\section{Summary}

In this work we have done a statistical analysis of some properties of
the contour lines of the Riemann zeta function.  We have uncovered some
amazing and previously unknown facets of the behavior of this remarkable
function.

\begin{acknowledgments}
The author thanks Stephen Wolfram for making the Wolfram CDF Player 8
available as a free public download.  He also thanks Peter Sarnak for
helpful conversations during the early stages of this work.  The
approximation of Eqn. (5) was pointed out to the author by an anonymous
referee.

\end{acknowledgments}



\end{document}